\documentclass[runningheads,a4paper]{llncs}

\usepackage{amssymb}
\usepackage{amsmath}

\usepackage{graphicx}
\usepackage{comment}

\usepackage{algorithmic}
\usepackage{algorithm} 
\usepackage{setspace}

\usepackage{listing}
\usepackage{listings}

\usepackage{xcolor}
\definecolor{mygreen}{RGB}{28,220,28} 
\definecolor{mylilas}{RGB}{170,55,241}

\usepackage{url}
\usepackage{hyperref}

\usepackage{tikz}
\usetikzlibrary{external}
\usepackage{tabularx}
\usepackage{tabulary}

\urldef{\mailsb}\path|jan.valdman@utia.cas.cz |
\urldef{\mailam}\path|alexmos@kma.zcu.cz |
\urldef{\mailcm}\path|matonoha@cs.cas.cz |

\newcommand{\nt}{|\mathcal{T}|}
\newcommand{\nn}{|\mathcal{N}|}

\newcommand{\R}{\mathbb R}

\def\IntO{\int\limits_{\Omega}}
\newcommand\dx{\, \mathrm{d}x}
\newcommand\dxy{\, \mathrm{d}x \mathrm{d}y}

\begin{document}

\mainmatter  
\title{Minimization of p-Laplacian via the Finite Element Method in MATLAB}

\titlerunning{Minimization of p-Laplacian via FEM in MATLAB}

\author{Ctirad Matonoha\inst{1}
\and
Alexej Moskovka\inst{2}
\and
Jan Valdman\inst{3,4}
\thanks{The work of C. Matonoha was supported by the long-term strategic development financing of the Institute of Computer Science (RVO:67985807). The work of A. Moskovka and J. Valdman was supported by the Czech-Austrian Mobility MSMT Grant: 8J21AT001.}
}

\authorrunning{Ctirad Matonoha, Alexej Moskovka, Jan Valdman}


\institute{$^1$The Czech Academy of Sciences, Institute of Computer Science, \\
Pod Vod\'{a}renskou v\v{e}\v{z}\'{\i}~2, 18207~Prague~8, Czechia \\
\mailcm \\
\vspace{0.15cm}
$^2$Faculty of Applied Sciences, Department of Mathematics, \\
University of West Bohemia, 
Technick\' a 8, 30614 Pilsen, Czechia \\
\mailam \\
\vspace{0.15cm}
$^3$The Czech Academy of Sciences, Institute of Information Theory \\ and Automation, 
Pod Vod\'{a}renskou v\v{e}\v{z}\'{\i}~4, 18208~Prague~8  \\
$^4$Department of Applied Informatics, Faculty of Science, University \\ of South Bohemia, Brani\v sovsk\' a 31, 37005~\v{C}esk\'{e}~Bud\v{e}jovice, Czechia 
\mailsb 
}
\maketitle
\begin{abstract}
Minimization of energy functionals is based on a discretization by the finite element method and optimization by the trust-region method. A key tool to an efficient implementation is a local evaluation of the approximated gradients together with sparsity of the resulting Hessian matrix. Vectorization concepts are explained for the p-Laplace problem in one and two space-dimensions.

\end{abstract}

\section{Introduction}
We are interested in a (weak) solution of the p-Laplace equation \cite{Drabek,Lindqvist}: 
\begin{equation} \label{pLapl}
\begin{split}
\Delta_p u &= f \qquad\quad \mbox{in} \:\: \Omega \, , \\
u &= g \qquad \,\, \mbox{on} \:\: \partial \Omega,
\end{split}
\end{equation}
where the p-Laplace operator is defined as $
\Delta_p u = \nabla \cdot \big( |\nabla u|^{p-2} \nabla u \big)$ for some power $p>1$. The domain $\Omega \in \mathbb{R}^d$ is assumed to have a Lipschitz boundary $\partial \Omega$, $f \in L^2(\Omega)$ and $g \in W^{1-1/p,p}(\partial \Omega)$, where $L$ and $W$ denote standard Lebesque and Sobolev spaces. It is known that \eqref{pLapl} represents an Euler-Lagrange equation corresponding to a minimization problem
\begin{equation} \label{energy1D}
J(u)=\min_{v \in V}J(v), \qquad J(v):=\frac{1}{p} \IntO |\nabla v|^p \dx - \IntO  f \, v  \dx,
\end{equation}
where $V=W^{1,p}_g(\Omega)=\{v \in W^{1,p}, v=g \mbox{ on } \partial \Omega \}$ includes Dirichlet boundary conditions on $\partial \Omega$. The minimizer $u \in V$ of \eqref{energy1D} is known to be unique for $p>1$. 
 
Due to the high complexity of the p-Laplace operator (with the exception of the case $p=2$ which corresponds to the classical Laplace operator), the analytical handling of \eqref{pLapl} is difficult.  The finite element method \cite{BarrettLiu,Ciarlet-FEM} can be applied as an approximation of  \eqref{energy1D} and results in a minimization problem
\begin{equation}\label{minimization_problem_discrete}
J(u_h)=\min_{v \in V_h} J(v), \qquad J(v):=\frac{1}{p} \IntO |\nabla v|^p \dx - \IntO  f \, v  \dx
\end{equation}
formulated over the finite-dimensional subspace $V_h$ of $V$. We consider for simplicity the case $V_h = P^1(\mathcal{T})$ only, where $P^1(\mathcal{T})$ is the space of  nodal basis functions defined on a triangulation $\mathcal{T}$ of the domain $\Omega$ using the simplest possible elements (intervals for $d=1$, triangles for $d=2$, tetrahedra for $d=3$). The subspace $V_h$ is spanned by a set of $n_b$ basis functions $\varphi_i(x) \in V_h, i=1, \dots, n_b$ and a trial function $v \in V_h$ is expressed by a linear combination 
$$v(x) = \sum_{i=1}^{n_b} v_i \, \varphi_i(x), \qquad x \in \Omega,$$
where $\bar v = (v_1, \dots, v_{n_b}) \in \mathbb{R}^{n_b}$ is a vector of coefficients. The minimizer $u_h \in V_h$ of \eqref{minimization_problem_discrete} is represented by a vector of coefficients $\bar u = (u_1, \dots, u_{n_b}) \in \mathbb{R}^{n_b}$ and some coefficients of $\bar u, \bar v$ related to Dirichlet boundary conditions are prescribed. 

In this paper, the first-order optimization methods are combined with FEM implementations \cite{AnjamValdman2015,RahmanValdman2013} in order to solve \eqref{minimization_problem_discrete} efficiently. These are the quasi-Newton (QN) and the trust-region (TR) methods \cite{conn2000} that are available in the MATLAB Optimization Toolbox.  The QN methods only require the knowledge of $J(v)$ and is therefore easily applicable. The TR methods additionally require the numerical gradient vector
$$\nabla J(\bar v) \in \mathbb{R}^{n_b}, \quad \bar v \in \mathbb{R}^{n_b}$$ and also allow to specify a sparsity pattern of the Hessian matrix 
$\nabla^2 J(\bar v) \in \mathbb{R}^{n_b \times n_b}, \bar v \in \mathbb{R}^{n_b}$, i.e., only positions (indices) of nonzero entries. The sparsity pattern is directly given by a finite element discretization. \\

We compare four different options:
\begin{itemize}
    \item[$\bullet$] option 1 : the TR method with the gradient evaluated directly via its explicit form and the specified Hessian sparsity pattern.
    \item[$\bullet$] option 2 : the TR method with the gradient evaluated approximately via central differences and the specified Hessian sparsity pattern.
    \item[$\bullet$] option 3 : the TR method with the gradient evaluated approximately via central differences and no Hessian sparsity pattern.
    \item[$\bullet$] option 4 : the QN method.
\end{itemize}

Clearly, option 1 is only applicable if the exact form of gradient is known while option 2 with the approximate gradient is only bounded to finite elements discretization and is always feasible.  Similarly to option 2, option 3 also operates with the approximate form of gradient, however the Hessian matrix is not specified. Option 3 serves as an intermediate step between options 2 and 4. Option 4 is based on the Broyden–Fletcher–Goldfarb–Shanno (BFGS) formula.

\section{One-dimensional problem}
The p-Laplace equation \eqref{pLapl} can be simplified as
\begin{equation} \label{pLapl1D}
\big( |u_x|^{p-2} u_x \big)_x = f \, \qquad \mbox{in} \:\: \Omega=(a,b)
\end{equation}
and the energy as 
\begin{equation} 
J(v):=\frac{1}{p} \int_{a}^b |v_x|^p \dx - \IntO  f \, v  \dx \, .
\end{equation}
\begin{figure}
\vspace{-0.8cm}
\centering
\includegraphics[width=0.71\textwidth]{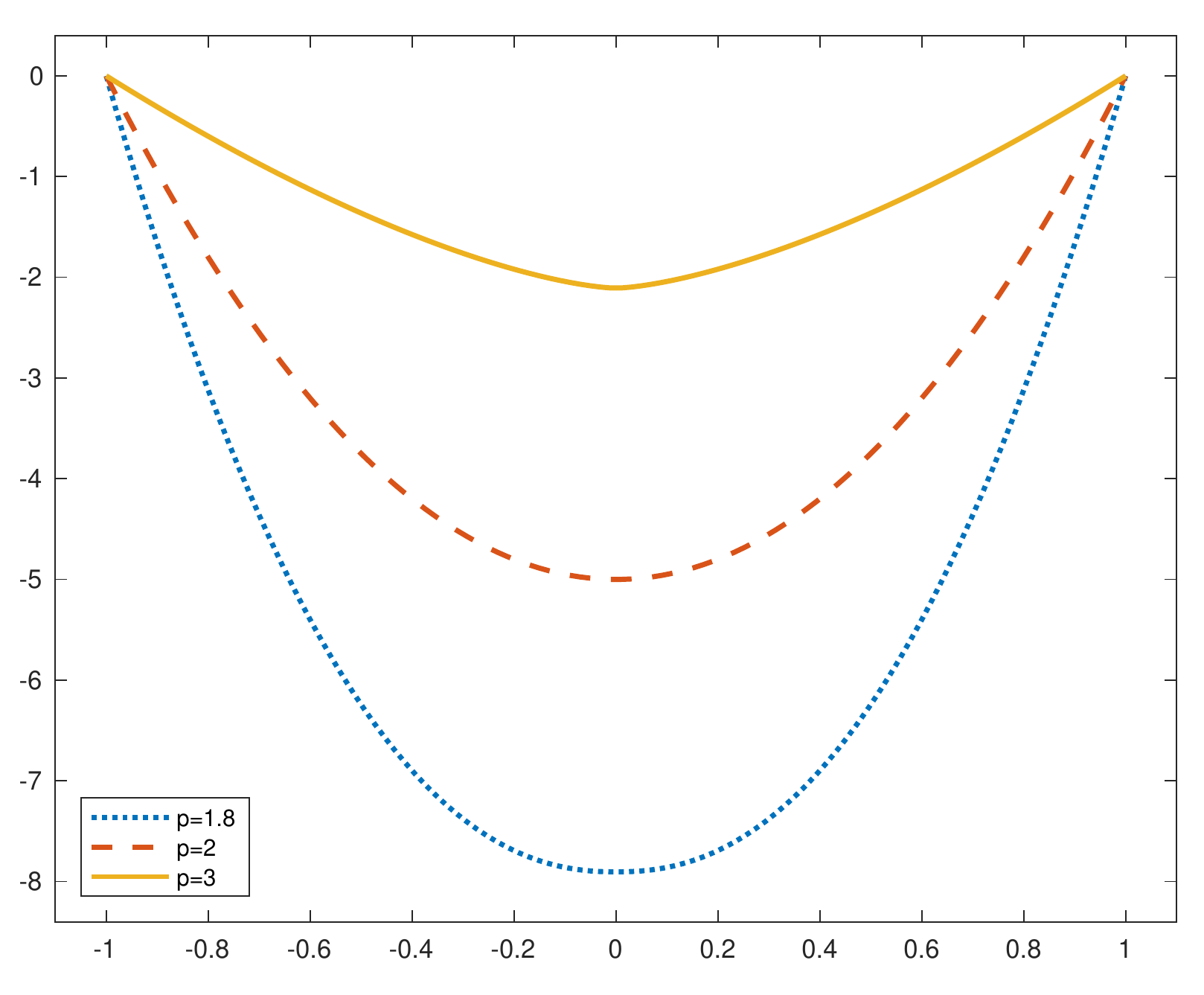}
  \caption{Solutions for $p \in \{1.8, 2, 3\}, \Omega=(-1,1), f=-10$ and Dirichlet boundary conditions $u(-1)=u(1)=0$.}
\label{pLaplace_1D}
\vspace{-0.5cm}
\end{figure}

Assume for simplicity an equidistant distribution of $n+2$ discretization points ordered in a vector
$(x_0, \ldots, x_{n+1}) \in \mathbb{R}^{n+2},$
where
$x_i:=a+i \, h$ for $i=0, 1, \ldots, n+1$ and 
$h:=(b-a)/(n+1)$
denotes an uniform length of all sub-intervals. It means that $x_0=a, x_{n+1}=b$ are boundary nodes. 

There are $n+2=n_b$ well-known hat basis functions 
$\varphi_0(x), \ldots, \varphi_{n+1}(x)$
satisfying the property 
$\varphi_i(x_j)= \delta_{ij}, i,j=0,\dots, n+1,$
where $\delta$ denotes the Kronecker symbol. 

Then, $v \in V_h$ 
is a piecewise linear and globally continuous function on $(a,b)$
represented by a vector of coefficients
$ \bar v=(v_0, \ldots, v_{n+1}) \in \mathbb{R}^{n+2}.$  
The minimizer $u_h \in V_h$ is similarly represented by a vector $\bar u=(u_0, \ldots, u_{n+1}) \in \mathbb{R}^{n+2}$. 
Dirichlet boundary conditions formulated at both interval ends imply 
$v_0=u_0=g(a), v_{n+1}=u_{n+1}=g(b),$ where boundary values $g(a), g(b)$ are prescribed.  

It is convenient to form a mass matrix $M \in \mathbb{R}^{(n+2) \times (n+2)}$ with entries
\begin{equation} \label{mass}
M_{i,j} = \int_{a}^b  \varphi_{i-1}(x) \, \varphi_{j-1}(x) \dx = h \cdot 
\begin{cases}
1/3, &  i = j \in \{1, n+2 \}  \\
2/3, &  i = j \in \{2,\dots, n+1 \} \\
1/6, & \mbox{$|i-j| = 1$} \\
0, & \mbox{otherwise}
\end{cases} .
\end{equation}
If we assume that $f \in V_h$ is represented by a (column) vector 
$ \bar f= (f_0, \ldots, f_{n+1}) \in \mathbb{R}^{n+2}$,
then the linear energy term reads exactly 
$$\int_a^b  f v  \dx = 
\bar f^T M \bar v = \bar b^T \bar v =  
\sum_{i=0}^{n+1} b_i v_i,$$
where $\bar b = (b_0, \ldots, b_{n+1}) = \bar{f}^T M \in \mathbb{R}^{n+2}$. 

The gradient energy term is based on the derivative $v_x$ which is a piecewise constant function and reads 
$$ v_x |_{(x_{i-1}, x_i)} = (v_{i}-v_{i-1})/{h}, \qquad i=1,\dots,n+1.
$$

 Now, it is easy to derive the following minimization problem:

\begin{problem}[p-Laplacian in 1D with Dirichlet conditions at both ends]\label{problem1} \\
Find $u = (u_1, \ldots, u_n) \in \mathbb{R}^n$ satisfying
\begin{equation}\label{pLaplace1Ddicrete} 
J(u) = \min_{v \in \mathbb{R}^n} J(v), \qquad J(v)=\frac{1}{p \,h^{p-1}} \sum_{i=1}^{n+1} |v_i - v_{i-1}|^p - \sum_{i=0}^{n+1} b_i v_i,
\end{equation}
where values $v_0:=g(a), v_{n+1}:=g(b)$ are prescribed. 
\end{problem}

Note that the full solution vector reads
$\bar u = (g(a), u, g(b)) \in \mathbb{R}^{n+2},$ where $u \in \mathbb{R}^{n}$ solves Problem \ref{problem1} above. 

Figure \ref{pLaplace_1D} illustrates discrete minimizers $\bar u$ for $(a, b)=(-1, 1)$, $f=-10$ and $p \in \{1.8, 2, 3\}$ assuming zero Dirichlet conditions $u(a)=u(b)=0$. Recall that the exact solution $u$ is known in this simple example. 

Table \ref{tab:laplace1D} depicts performance of all four options for the case $p=3$ only, in which the exact energy reads $J(u) = -\frac{16}{3}\sqrt{10} \approx -16.8655$. 
The first column of every option shows evaluation time, while the second column provides the total number of linear systems to be solved (iterations), including rejected steps.
Clearly, performance of options 1 and 2 dominates over options 3 and 4.  
\begin{table}[H]
    \centering
    \begin{tabularx}{0.99\textwidth}
    {r
    |>{\raggedleft\arraybackslash}X 
    >{\raggedleft\arraybackslash}X 
    |>{\raggedleft\arraybackslash}X
    >{\raggedleft\arraybackslash}X
    |>{\raggedleft\arraybackslash}X 
    >{\raggedleft\arraybackslash}X 
    |>{\raggedleft\arraybackslash}X
    >{\raggedleft\arraybackslash}X
    }
      &  \multicolumn{2}{c|}{option 1:} & \multicolumn{2}{c|}{option 2:} & \multicolumn{2}{c|}{option 3:} & \multicolumn{2}{c}{option 4:}   \\
      \hline
     n$\,$   &  time & iters & time & iters & time  & iters & time & iters \\
     \hline
 1e1 &   0.01 &    8  &   0.01 &    6  &   0.02 &    6  &   0.02 &   17  \\ 
 1e2 &   0.03 &   12  &   0.05 &   11  &   0.49 &   11  &   0.29 &   94  \\ 
 1e3 &   0.47 &   37  &   0.50 &   15  &  96.22 &   14  &  70.51 &  922  \\ 
    \end{tabularx}
     \vspace{0.2cm}
    \caption{
    MATLAB performance in 1D for $p=3$. Times are given in seconds.
    }\label{tab:laplace1D}
    \vspace{-0.5cm}
\end{table}

\section{Two-dimensional problem}
The equation \eqref{pLapl} in 2D has the form
\begin{equation}
\nabla \cdot \Bigg( \bigg[ \Big(\frac{\partial u}{\partial x}\Big)^2 + \Big(\frac{\partial u}{\partial y}\Big)^2 \bigg]^{\frac{p-2}{2}} \nabla u \Bigg) = f \qquad \mbox{in} \:\: \Omega
\end{equation}
and the corresponding energy reads
\begin{equation} 
J(v):=\frac{1}{p} \iint\limits_{\Omega} \Big(|v_x|^p + |v_y|^p \Big) \dxy - \iint\limits_{\Omega}  f \, v  \dxy \, .
\end{equation}
\begin{figure}
\vspace{-0.8cm}
\centering
\begin{minipage}[c]{0.48\textwidth}
\includegraphics[width=\textwidth]{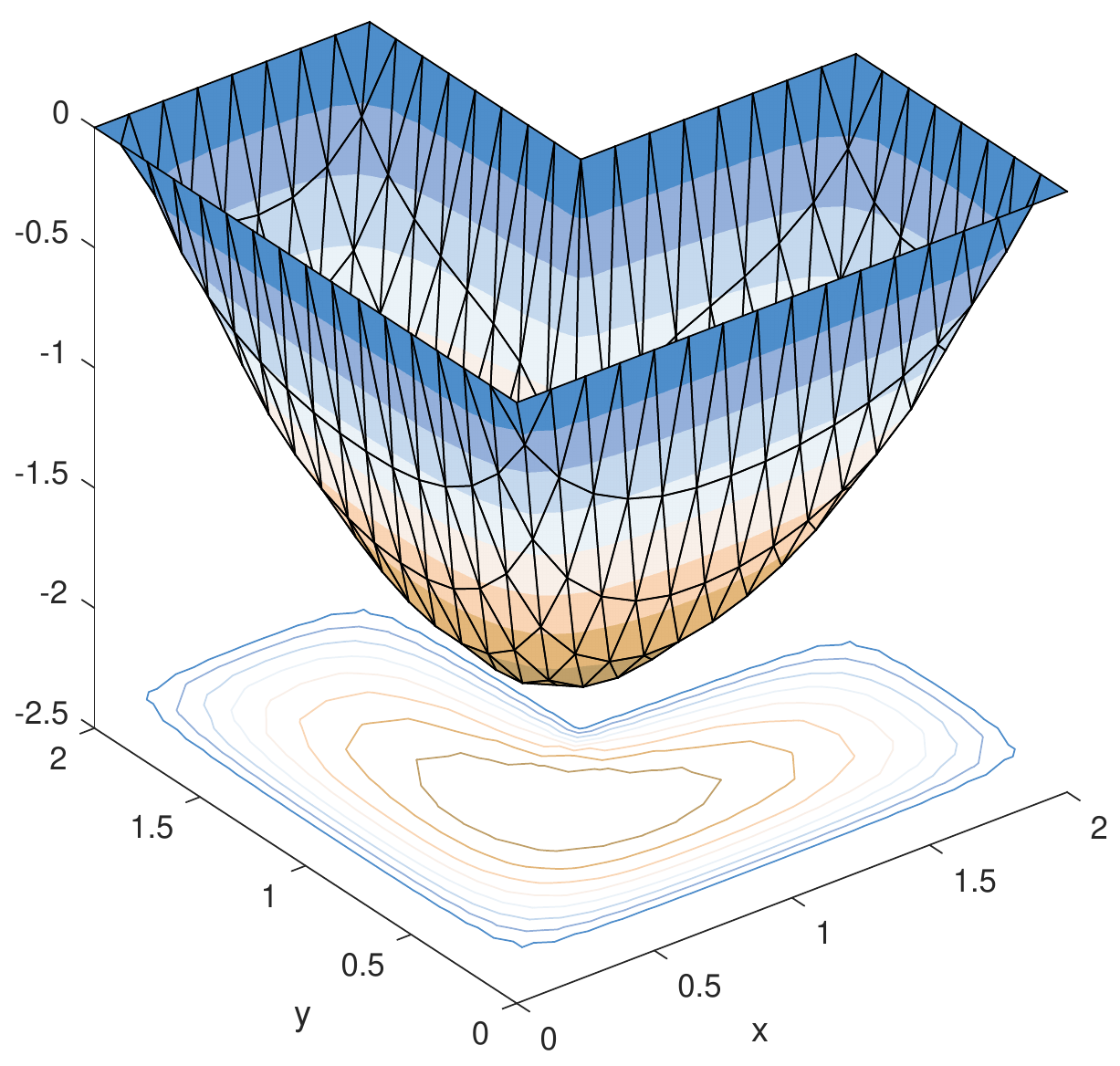}
\end{minipage} \:\:\:
\vspace{0.01\textwidth}
\begin{minipage}[c]{0.48\textwidth}
\includegraphics[width=\textwidth]{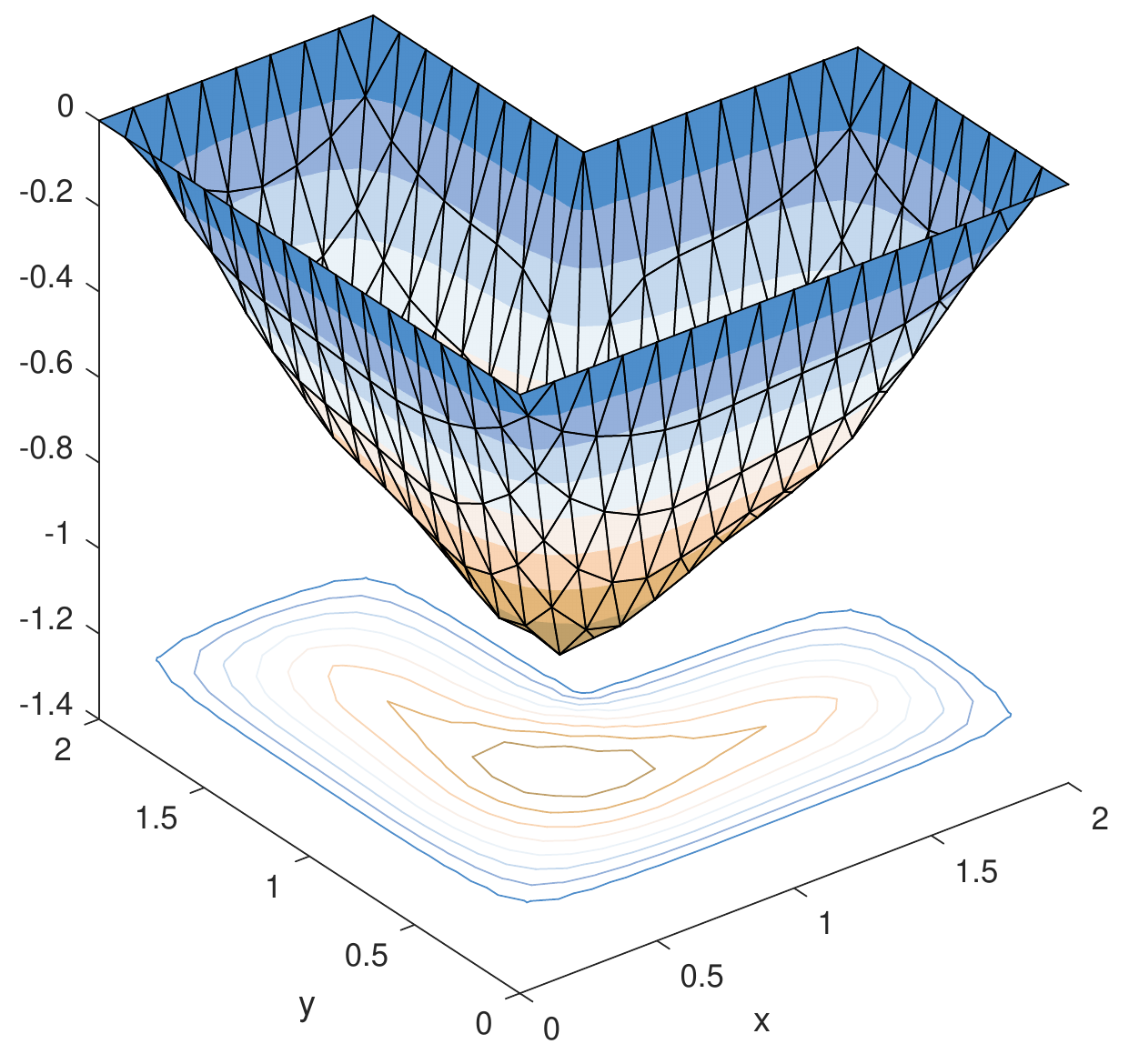}
\end{minipage}
\caption{
Numerical solutions with contour lines for $p=1.8$ (left) and $p=3$ (right) and a L-shape domain $\Omega, f=-10$ and zero Dirichlet boundary conditions on $\partial \Omega$.}
\label{pLaplace_2D}
\vspace{-0.5cm}
\end{figure}

Assume a domain $\Omega \in \mathbb{R}^2$ with a polygonal boundary $\partial \Omega$ is discretized by a regular triangulation of triangles \cite{Ciarlet-FEM}. The sets $\mathcal{T}$ and $\mathcal{N}$ denote the sets of all triangles and their nodes (vertices) and $\nt$ and $\nn$ their sizes, respectively. 
Let \( \mathcal{N}_{dof} \subset \mathcal{N} \) be the set of all internal nodes and \( \mathcal{N} \, \backslash \, \mathcal{N}_{dof} \) denotes the set of boundary nodes. 

A trial function $v \in V_h = P_1(\mathcal{T})$ is a globally continuous and linear scalar function on each triangle $T \in \mathcal{T}$ represented by a vector of coefficients 
$ \bar v = (v_1, \ldots, v_{\nn}) \in \mathbb{R}^{\nn}.$
Similarly the minimizer $u_h \in V_h$ is represented by a vector of coefficients
$\bar u = (u_1, \ldots, u_{\nn}) \in \mathbb{R}^{\nn}.$
Dirichlet boundary conditions imply
\begin{equation}\label{dirichlet}
v_i=u_i=g(N_i) \, , \qquad \mbox{ where } N_i \in \mathcal{N} \, \backslash \, \mathcal{N}_{dof},    
\end{equation}
and the function $g: \partial \Omega \rightarrow \mathbb{R}$ prescribes 
Dirichlet boundary values.
\begin{example}
A triangulation $\mathcal{T}$ of the L-shape domain $\Omega$ is given in Figure \ref{mesh_sparsity} (left) in which $\nt=24, \nn=21.$ The Hessian sparsity pattern (right) can be directly extracted from the triangulation: it has a nonzero value at the position $i,j$, if nodes $i$ and $j$ share a common edge. 
\begin{figure}[h]
\vspace{-0.2cm}
\centering
\begin{minipage}[c]{0.48\textwidth}
\includegraphics[width=0.99\textwidth]{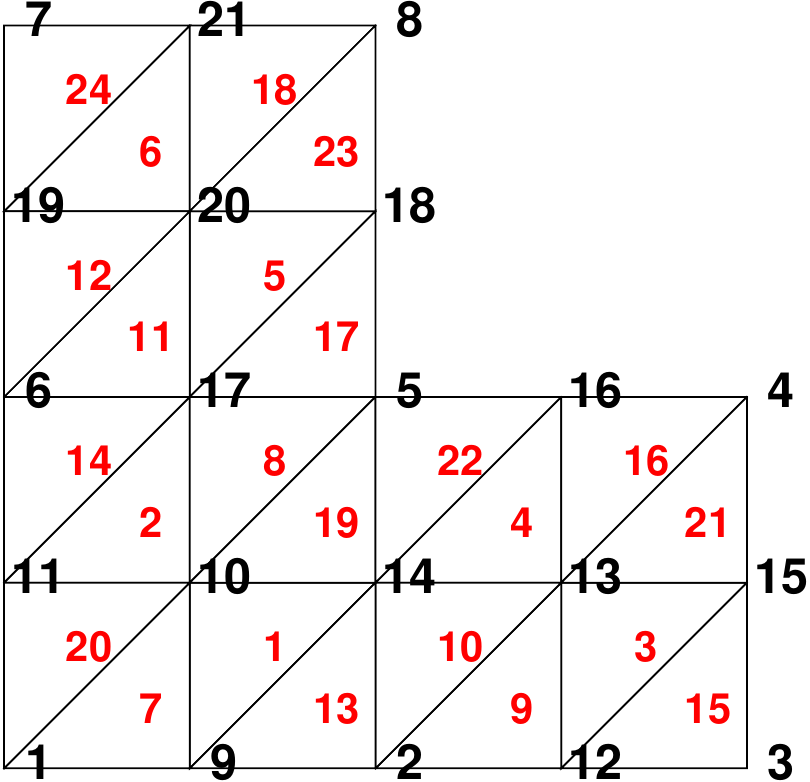}
\end{minipage}
\begin{minipage}[c]{0.48\textwidth}
\includegraphics[width=0.99\textwidth]{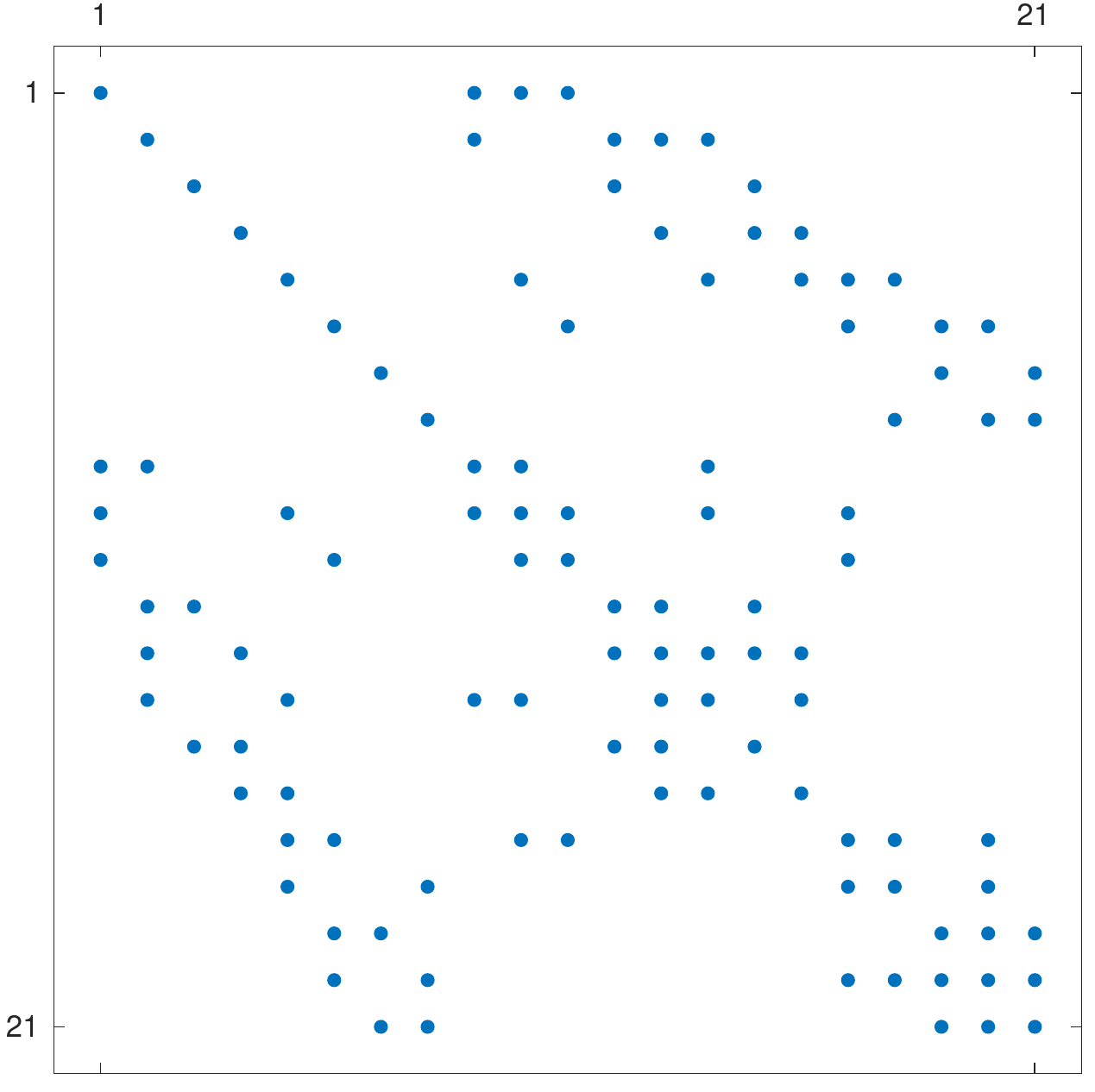}
\end{minipage}
\caption{A triangular mesh (left) and the corresponding Hessian sparsity pattern (right).}
\label{mesh_sparsity}
\vspace{-0.3cm}
\end{figure}
The set of internal nodes that appear in the minimization process reads 
$\mathcal{N}_{dof} = \{N_{10},N_{13}, N_{14}, N_{17},N_{20} \},$
while the remaining nodes belong to the boundary $\partial \Omega$. 

\end{example}

For an arbitrary node $N_k$, $k \in \{1, 2, \hdots, \nn\}$ we define a global basis function $\varphi_k$ which is linear on every triangle and holds
$\varphi_k(N_l) = \delta_{k l}, l \in \{1, 2, \hdots, \nn\} \, .
$
Note that with these properties all global basis functions are uniquely defined.

Similarly to 1D, assume $f \in V_h$ is represented by a (column) vector
$ \bar{f} \in \mathbb{R}^{\nn}$, and introduce a (symmetric)
mass matrix $M \in \mathbb{R}^{\nn \times \nn}$ with entries
$M_{i,j}  = \iint\limits_{\Omega} \varphi_i \varphi_j \dxy.$
Then it holds $\iint\limits_{\Omega}  f v  \dxy = \sum_{i=1}^{\nn} b_i v_i,$ where $b = \bar{f}^T M \in \mathbb{R}^{\nn}$.

Next, for an arbitrary element $T_i \in \mathcal{T}$, $i \in \{1,2,\hdots, \nt\}$, denote 
$\varphi^{i,1}, \varphi^{i,2}, \varphi^{i,3}$ all three local basis functions on the $i$-th element and let
$    \varphi_x^{i,j}, \varphi_y^{i,j}, j \in \{1,2,3\}
$
be the partial derivatives with respect to 'x' and 'y' of the $j$-th local basis function on the $i$-th element, respectively.
In order to formulate the counterpart of  \eqref{pLaplace1Ddicrete} in two dimensions, we define gradient vectors $v_{x,el}, v_{y,el} \in \R^{\nt}$ with entries
\begin{equation*}
v_{x,el}^{i} = \sum_{j=1}^3 \varphi_{x}^{i,j} v^{i,j}, \qquad 
v_{y,el}^{i} = \sum_{j=1}^3 \varphi_{y}^{i,j} v^{i,j} \, ,
\end{equation*}
where $v^{i,j}$ is the value of $v$ in the $j$-th node of the $i$-th element.

With these substitutions we derive the 2D counterpart of Problem \ref{problem1}:
\begin{problem}[p-Laplacian in 2D]\label{problem3}
Find a minimizer $u \in \mathbb{R}^{\nn}$ satisfying
\begin{equation}\label{pLaplace2Ddicrete}
J(u) = \min_{v \in \mathbb{R}^{\nn}} J(v), \quad J(v) = \frac{1}{p} \sum_{i=1}^{\nt} |T_i| \Big(
|v_{x,el}^i|^p + |v_{y,el}^i|^p \Big) - \sum_{i=1}^{\nn} b_i v_i
\end{equation}
with prescribed values $v_i = g(N_i)$ for $N_i \in \mathcal{N} \, \backslash \, \mathcal{N}_{dof}$. 
\end{problem}

\begin{table}[h]
 \vspace{-0.5cm}
    \centering
    \begin{tabularx}{0.99\textwidth}
    {>{\raggedleft\arraybackslash}X
    |>{\raggedleft\arraybackslash}X 
    >{\raggedleft\arraybackslash}X 
    |>{\raggedleft\arraybackslash}X
    >{\raggedleft\arraybackslash}X
    |>{\raggedleft\arraybackslash}X 
    >{\raggedleft\arraybackslash}X 
    |>{\raggedleft\arraybackslash}X
    >{\raggedleft\arraybackslash}X
    }
      &  \multicolumn{2}{c|}{option 1:} & \multicolumn{2}{c|}{option 2:} & \multicolumn{2}{c|}{option 3:} & \multicolumn{2}{c}{option 4:}   \\
      \hline
     \( |\mathcal{N}_{dof}| \)  &  time & iters & time & iters & time  & iters & time & iters \\
     \hline
 33 &   0.04 &    8  &   0.05 &    8  &   0.15 &    8  &   0.06 &   19  \\ 
 161 &   0.20 &   10  &   0.29 &    9  &   3.19 &    9  &   0.56 &   31  \\ 
 705 &   0.75 &    9  &   1.17 &    9  &  70.59 &    9  &  12.89 &   64  \\ 
 2945 &   3.30 &   10  &   5.02 &    9  &  - & -  &  388.26 &  133  \\ 
 12033 &  16.87 &   12  &  24.07 &   10  &  - & -  &  - & -  \\ 
 48641 &  75.32 &   12  &  107.38 &   10  &  - & -  &  - & -  \\ 
    \end{tabularx}
     \vspace{0.2cm}
    \caption{MATLAB performance in 2D for $p=3$. Times are given in seconds.
    }\label{tab:laplace2D}
    \vspace{-0.5cm}
\end{table}

Figure \ref{pLaplace_2D} illustrates numerical solutions for the L-shape domain from Figure \ref{mesh_sparsity}, for $f=-10$ and $p \in \{1.8, 3 \}$. 
Table \ref{tab:laplace2D} depicts performance of all options for $p=3$. Similarly to 1D case (cf. Table \ref{tab:laplace1D}), performance of options 1 and 2 clearly dominates over options 3 and 4. Symbol '-' denotes calculation which ran out of time or out of memory. The exact solution $u$ is not known in this example but numerical approximations provide the upper bound $J(u) \approx -8.1625$.

\subsection{Remarks on 2D implementation}
As an example of our MATLAB implementation, we introduce below the following block describing the evaluation of formula \eqref{pLaplace2Ddicrete}:
\begin{listing}
\begin{lstlisting}
function e=energy(v)
    v_elems=v(elems2nodes);                           
    v_x_elems=sum(dphi_x.*v_elems,2);                
    v_y_elems=sum(dphi_y.*v_elems,2);           
    intgrds=(1/p)*sum(abs([v_x_elems v_y_elems]).^p,2);
    e=sum(areas.*intgrds) - b'*v;
end
\end{lstlisting}
\end{listing}
The whole code is based on several matrices and vectors that contain the topology of the triangulation and gradients of basis functions. Note that these objects are assembled effectively by using vectorization techniques from \cite{AnjamValdman2015,RahmanValdman2013} once and do not change during the minimization process. These are (with dimensions): 

\begin{description}
    \itemsep 0.5mm
    \item[\texttt{elems2nodes}] $\nt \times 3$ - for a given element returns three corresponding nodes
    \item[\texttt{areas}] $\nt \times 1$ - vector of areas of all elements, $\texttt{areas(i)} = |T_i|$
    \item[\texttt{dphi\_x}] $\nt \times 3$ - partial derivatives of all three basis functions with respect to $x$ on every element
    \item[\texttt{dphi\_y}] $\nt \times 3$ - partial derivatives of all three basis functions with respect to $y$ on every element
\end{description}

\noindent The remaining objects are recomputed in every new evaluation of the energy:
\begin{description}
    \itemsep 0.5mm
    \item[\texttt{v\_elems}] $\nt \times 3$ - where \texttt{v\_elems(i,j)} represents $v^{i,j}$ above
    \item[\texttt{v\_x\_elems}] $\nt \times 1$ - where \texttt{v\_x\_elems(i)} represents $v_{x,el}^i$ above
    \item[\texttt{v\_y\_elems}] $\nt \times 1$ - where \texttt{v\_y\_elems(i)} represents $v_{y,el}^i$ above
\end{description}

The evaluation of the energy above is vital to option 4. For other options, exact and approximate gradients of the discrete energy \eqref{pLaplace2Ddicrete} are needed, but not explained in detail here. 
Additionally, for options 1 and 2, the Hessian pattern is needed and is directly extracted from the object \texttt{elems2nodes} introduced above.

\subsubsection*{Implementation and outlooks}
Our MATLAB implementation is available 
at
\begin{center}
\url{https://www.mathworks.com/matlabcentral/fileexchange/87944} 
\end{center}
for download and testing. The code is designed in a modular way that different scalar problems involving the first gradient energy terms can be easily added. Additional implementation details on evaluation of exact and approximate gradients will be explained in the forthcoming paper. 

We are particularly interested in further vectorization of current codes resulting in faster performance and also in extension to vector problems such as nonlinear elasticity. 
Another goal is to exploit line search methods from \cite{ufo2017}.



\bibliographystyle{abbrv}

\end{document}